\documentclass[12 pt]{article}
\usepackage{latexsym}
\usepackage{amssymb}
\usepackage{amsthm} 
\usepackage{amsmath}
\usepackage[]{color}
\usepackage{enumitem}
\usepackage{hyperref}
\usepackage{authblk}

\usepackage{tikz}

\newtheorem{thm}{Theorem}

\newtheorem{corollary}{Corollary}

\newtheorem{lemma}{Lemma}

\theoremstyle{definition}
\newtheorem*{defn}{Definition}
\newtheorem{example}{Example}

\newtheorem{problem}{Problem}

\newcommand\dist{{\rm Dist}}
\renewcommand\det{{\rm Det}}
\newcommand\aut{{\rm Aut}}
\newcommand\setstab{{\rm SetStab}}

\newcommand\fdist{{\rm Fdist}}
\newcommand\lrho{{\rho^\ell}}
\newcommand\urho{{\rho^u}}
\newcommand\1{\bigskip\noindent}
\renewcommand\th{^{th}}

\newcommand\st{^{st}}

\newcommand\symm{{\rm Symm}}

\title{Paint Cost and the Frugal Distinguishing Number}

\author{\Large Debra L. Boutin,  \\ \footnotesize Hamilton College, Clinton, NY  13323\\ {\tt dboutin@hamilton.edu}}

\date{\today}

\begin{document}

\maketitle

\begin{abstract}  You are handed a graph with vertices in a neutral color and asked to color a subset of vertices with expensive paints in $d$ colors in such a way that only the trivial symmetry preserves the color classes.  Your goal is to minimize the number of vertices needing this expensive paint.  This paper address the issues surrounding your choices.  In particular, a graph is said to be \emph{$d${-}distinguishable} if there exists a coloring with $d$ colors so that only the trivial automorphism preserves the color classes.  The \emph{distinguishing number} of $G$, denoted $\dist(G)$, is the smallest $d$ for which $G$ is $d${-}distinguishable. We define the \emph{paint cost of $d${-}distinguishing}, denoted $\rho^d(G)$,  to be  the minimum number of vertices that need to be painted to $d${-}distinguish $G$.  This cost varies with $d$.  The maximum paint cost for $G$ is called the \emph{upper paint cost}, denoted $\urho(G)$, and occurs when $d{=}\dist(G)$; the minimum paint cost is called the \emph{lower paint cost}, denoted $\rho^\ell(G)$.  Further, we define the smallest $d$ for which the paint cost is $\rho^\ell(G)$, to be the  \emph{frugal distinguishing number} $\fdist(G)$.  In this paper we formally define $\rho^d(G)$, $\urho(G)$, $\lrho(G)$, and $\fdist(G)$.  We also show that $\urho(G)$ and $\lrho(G)$, as well as $\fdist(G)$ and $\dist(G)$, can be arbitrarily large multiples of each other.  Lastly,  we find these parameters for the book graph $B_{m,n}$, summarized as follows.  For $n\geq 2$ and  $m\geq 4$, we show
\begin{itemize}

\item $\lrho(B_{m,n}) =  n{-}1;$

\item $\rho^u(B_{m,n}) \geq (m{-}2) \left( n-k^{m-3} \right) +1, \text{ where }  k=\dist(B_{m,n});$

\item $\fdist(B_{m,n}) = 2+\left\lfloor \frac{n{-}1}{m{-}2} \right\rfloor.$

\end{itemize}  \end{abstract}

\1{\bf Keywords:} graph distinguishing, cost of distinguishing, determining number, book graph

\section{Introduction}

Graph distinguishing is the act of coloring the vertices of a graph in such a way that only the trivial automorphism preserves the color classes. Some call this process \emph{symmetry breaking}. More specifically, a coloring of the vertices of a graph $G$ with  $d$ colors is called a \emph{$d${-}distinguishing coloring} if only the trivial automorphism preserves the color classes.   A graph is called \emph{$d${-}distinguishable} if it has a $d${-}distinguishing coloring. The \emph{distinguishing number} of $G$, denoted $\dist(G)$, is the smallest $d$ necessary for a $d$-distinguishing coloring of $G$.  This concept has also been studied under the name \emph{asymmetric coloring}.\medskip

There has been substantial interest in graph distinguishing in the last few decades. Much of the work proves that for a large number of graph families all but a finite number of members are 2{-}distinguishable. Examples of such families of finite graphs include: hypercubes $Q_n$ with $n\geq 4$~\cite{BoCo2004}, Cartesian powers $G^n$ for a connected graph $G\ne K_2,K_3$ and $n\geq 2$~\cite{Al2005, ImKl2006,KlZh2007}, and Kneser graphs $K_{n:k}$ with $n\geq 6, k\geq 2$~\cite{AlBo2007}. Examples of such families of infinite graphs include: the denumerable random graph~\cite{ImKlTr2007}, the infinite hypercube~\cite{ImKlTr2007}, locally finite trees with no vertex of degree 1 \cite{WaZh2007}, and denumerable vertex{-}transitive graphs of connectivity 1~\cite{SmTuWa2012}. \medskip

In 2008, Boutin~\cite{Bo2008} defined the \emph{cost of 2{-}distinguishing} a 2{-}distinguishable graph $G$, denoted $\rho(G)$,  to be the minimum size of a color class over all 2{-}distinguishing colorings of $G$.   There has been moderate and increasing interest in studying the cost of 2{-}distinguishing in the past few years. Graph families with known or bounded 2-distinguishing cost are: hypercubes \cite{Bo2020}, Kneser graphs  \cite{Bo2013b}, and Split Praeger-Xu graphs \cite{ILTW2020}.\medskip

 Note that the cost of 2-distinguishing measures both how close $G$ is to being 1-distinguishable, as well as the number of vertices that need (re)coloring to achieve a $2$-distinguishing coloring.  In generalizing the concept of cost  to more general $d$-distinguishing colorings, we find that we must separate these two measurement goals.  We can ask how close $G$ is to being $(d{-}1)$-distinguishable (which is the size of a smallest color class) or we can ask how many vertices we need to (re)color to achieve a $d$-distinguishing coloring (which is the size of the smallest complement of a color class). In~\cite{AlSo2021a}, Alikhani and Soltani take the first approach by defining the \emph{cost number}, denoted $\rho_d(G)$, to be the size of a smallest color class over all $d$-distinguishing colorings of $G$.  In this paper we take the second approach by defining the \emph{paint cost of $d${-}distinguishing} $G$, or equivalently the \emph{$d$-paint cost}, denoted $\rho^d(G)$, to be the minimum size of the complement of a color class over all $d$-distinguishing colorings.  As one expects from generalizations, if $G$ is 2-distinguishable then $\rho_2(G)=\rho(G)=\rho^2(G)$.\medskip

\begin{example} The cycle $C_5$ has distinguishing number 3.  One 3-distinguishing coloring of $C_5$ is, in cyclic order,  red, green, grey, grey, grey.  It is easy to argue that $\rho_3(C_5) =1$ and that $\rho^3(C_5)=2$.  That is, $C_5$ is one vertex away from being $2$-distinguishable, but it takes coloring at least two vertices to achieve a 3-distinguishing coloring.\end{example}

Determining sets are useful in finding the distinguishing number and the paint cost of $d$-distinguishing. A subset $S\subseteq V(G)$ is said to be a \emph{determining set} for $G$ if  only the trivial automorphism fixes its elements pointwise.  The \emph{determining number} of a graph $G$, denoted $\det(G)$, is the minimum size of a determining set. A determining set of minimum size is called a \emph{minimum determining set}.  Intuitively, if we think of automorphisms of a graph as allowing vertices to move among themselves, one can think of the determining number as the smallest number of vertices that need to be \lq\lq pinned down" to fix the entire graph.   The determining number has also been studied under the names \emph{fixing number} and \emph{rigidity index}.\medskip

Though distinguishing numbers and determining numbers were introduced by different people and for different purposes, they have strong connections. In~\cite{AlBo2007}, Albertson and Boutin show that if $G$ has a determining set $S$ of size $d$, then coloring each vertex in $S$ a distinct color from $\{1,\ldots,d\}$, and coloring the remaining vertices with color $d{+}1$, yields a $(d{+}1)${-}distinguishing coloring of $G$. Thus, $\dist(G) \leq \det(G) {+}1$. Further, in~\cite{Bo2013a} Boutin proves that in a 2{-}distinguishing coloring of $G$, each color class is a determining set, though not necessarily of minimum size. Thus, if $G$ is 2{-}distinguishable, then $\det(G) \leq \rho(G)$.\medskip

Notice that in defining $d$-paint cost, we do not explicitly use $\dist(G)$.  We only use the fact that $G$ is $d${-}distinguishable.  Thus  $d$-paint cost is a valid parameter for any $d\geq \dist(G)$. Further, $\rho^d(G)$ varies as $d$ does.  This raises one of the main points of this paper.  What is the largest value for $\rho^d(G)$? Call this the \emph{upper paint cost} and denote it by $\rho^u(G)$.  In Section~\ref{sec:tools}, we see that for all $d\geq \dist(G)$, we have $\det(G)\leq \rho^d(G)$, and further that there exist $d$ for which equality holds.  Thus, another question is, what is the smallest $d$ for which $\rho^d(G)=\det(G)$?  Call this the \emph{frugal distinguishing number}, and denote it by $\fdist(G)$.\medskip

This paper is organized as follows.  In Section~\ref{sec:tools}, we more formally define the various types of distinguishing costs, introduce frugal graph distinguishing, and provide a few examples.  In Section~\ref{sec:books} we study these symmetry parameters on the book graph $B_{m,n}$, which consist of $n$ copies of $C_m$ identified along an edge.  There we see that, for $n\geq 2$ and $m\geq 4$, 

\begin{itemize}

\item $\det(B_{m,n}) =  n{-}1;$

\item $\rho^u(B_{m,n}) \geq (m{-}2) \left( n-k^{m-3} \right) +1, \text{ where }  k=\dist(B_{m,n});$

\item $\fdist(B_{m,n}) = 2+\left\lfloor \frac{n{-}1}{m{-}2} \right\rfloor.$

\end{itemize}\medskip

\section{Definitions, Tools, and Examples}\label{sec:tools}

All graphs in this paper are finite and simple. Given a graph $G$ and $S~\subseteq~V(G)$, denote the complement of $S$ in $V(G)$ by $S^c$. Denote by $\symm(X)$ the set of permutations of the elements of set $X$.  We now begin more formally defining some of our fundamental terms. 

\begin{defn}  Let $G$ be a $d${-}distinguishable graph.  The \emph{paint cost} \emph{of} \emph{$d$-distinguishing} $G$, or equivalently the \emph{$d$-paint cost}, denoted $\rho^d(G)$, is the minimum size of the complement of a color class over all $d${-}distinguishing colorings for $G$.\end{defn}

Of course the minimum size of a complement of a color class occurs precisely when the size of the color class itself is  maximum; and the latter may be easier to compute.  Denote by $R^d(G)$ the maximum size of a color class over all $d${-}distinguishing colorings for $G$.  Then $\rho^d(G) = |V(G)|{-}R^d(G)$. 

\begin{lemma} \label{lem:disdet} \rm If $S$ is the complement of a color class in a distinguishing coloring for $G$, then $S$ is a determining set for $G$.\end{lemma}

\begin{proof}  Suppose we have a $d$-distinguishing coloring for $G$, and $T$ a color class of that coloring.  Denote $T^c$ by $S$.   As the complement of a color class, $S$ is the union of the remaining $d{-}1$ color classes in the distinguishing coloring. If $\alpha\in \aut(G)$ fixes $S$ pointwise, then $\alpha$ (trivially) preserves the $d{-}1$ color classes in $S$.  Thus $\alpha$ also preserves the $d\th$ color class $T$. Now, since $\alpha$ preserves all the color classes of a distinguishing coloring, $\alpha$ is the identity automorphism.  Thus $S$ is a determining set for $G$.\end{proof}

\begin{thm} \label{thm:detdist}\rm If $G$ is $d$-distinguishable, then $\det(G) \leq \rho^d(G)$, and equality can be achieved.  \end{thm}

\begin{proof} By definition of $\rho^d(G)$, there exists a $d$-distinguishing coloring of $G$ with a color class of size $|V(G)|-\rho^d(G)$.  Let $S$ be the complement of this color class.  By definition of $S$, $|S| = \rho^d(G)$.  Further, by Lemma~\ref{lem:disdet}, $S$ is a determining set for $G$.  Thus $\det(G)\leq \rho^d(G)$.\medskip

To show that equality is achievable, consider the following. Let $S$ be a minimum determining set for $G$.  Let $d=\det(G)$.  Color each vertex of $S$ distinctly with one of $d$ colors.  Color the vertices of $S^c$ with a $(d+1)\st$ color.  Since all vertices in $S$ have distinct colors, if $\alpha\in \aut(G)$ preserves the color classes, then it preserves $S$ pointwise.  Since $S$ is a determining set for $G$, this tells us that $\alpha$ is the trivial automorphism.  Thus the coloring is $(d{+}1)$-distinguishing with $(d{+}1)$-paint cost precisely $\det(G)$.  Since no $(\det(G){+}1)$-distinguishing coloring have paint cost less than $\det(G)$, we get that $\rho^{\det(G)+1}(G)=\det(G)$.\end{proof}

If we are given $d{+}1$ colors, we could choose to use only $d$ of them and achieve the same size color classes as with $d$ colors.  However, having an extra color may potentially allow us to color fewer vertices.  This yields the following fact.

\begin{thm} \rm If $d\geq \dist(G)$, then $\rho^d(G) \geq \rho^{d{+}1}(G)$.\end{thm}

 Given the inequality $\det(G)\leq \rho^d(G)$ and the fact that equality can be achieved, a natural question is: What is the smallest $d$ for which equality holds?  This motivates the following definition.

\begin{defn}  The \emph{frugal distinguishing number} of a graph $G$, denoted $\fdist(G)$, is the smallest $d$ for which $\rho^d(G) = \det(G)$. \end{defn}

Theorem~\ref{thm:detdist} and the definition above combine to give us the following corollary.

\begin{corollary} Let $D=\dist(G)$ and $F=\fdist(G)$.  Then $$\rho^D(G) \geq \rho^{D{+}1}(G) \geq \cdots \geq \rho^{F{-}1}(G) \geq \rho^F(G)= \rho^{F{+}1}(G) = \cdots = \det(G).$$  \end{corollary}

Since $\rho^{\dist(G)}(G)$ is the largest possible paint cost, we call it the \emph{upper paint cost}, denoted $\rho^u(G)$, and since $\rho^{\fdist(G)}(G)$ is the smallest possible paint cost, we call it the \emph{lower paint cost}, denoted $\rho^\ell(G)$.\medskip

One approach to finding $\fdist(G)$ is to begin with a minimum determining set $S$, and  give it a distinguishing coloring.  This is a more refined version  of what we did in the second half of the proof of Theorem~\ref{thm:detdist}.  However, we must first formalize what it means to distinguish a set of vertices.\medskip

Recall that for a subset $S\subseteq V(G)$, the \emph{set stabilizer} of $S$, denoted $\setstab(S)$, is the set of automorphisms of $G$ that preserve $S$ setwise. 

\begin{defn} A vertex coloring of $S\subseteq V(G)$ is said to be \emph{set distinguishing} if whenever $\alpha\in \setstab(S)$ preserves the color classes of $S$, then $\alpha$ fixes the elements of $S$ pointwise. \end{defn} 

It is straightforward to show that for any set $X$, the action of $\symm(X)$ on $X$ requires $|X|$ distinct colors to distinguish.   We use this fact a number of times in the remainder of the paper.  Also, note that if $\alpha\in \setstab(S)$ preserves a distinguishing coloring of $S$, it need not fix $S^c$ pointwise, only  preserve it setwise.  The following example illustrates this.

\begin{example} Consider $K_{5,1}$ with $V(K_{5,1})=\{0,1,2,3,4,5\}$ labeled so that vertex 0 has degree 5 and the remainder have degree 1. Let $S=\{1,2,3\}$.   It is easy to see that $\aut(K_{5,1}) =\symm(\{1,2,3,4,5\})$ and that $\setstab(S)= \symm(\{1,2,3\}) \times \symm(\{4,5\})$.  We can distinguish $S$ by giving each of its vertices a distinct color.  The nontrivial automorphism $\alpha\in \setstab(S)$ that permutes vertices 4 and 5  while fixing all other vertices, preserves the color classes of $S$, fixes the vertices of $S$ pointwise, and preserves $S^c$ setwise.\end{example}

As we apply the idea of distinguishing a minimum determining set, it is important to note that not all minimum determining sets have the same, or isomorphic, set stabilizers.  In particular, some minimum determining sets may require more colors to distinguish than others. Our goal is to find a minimum determining set with minimal symmetry and to break those symmetries with a minimum number of colors.  The following example demonstrates differences in  set stabilizers of minimum determining sets and the effect on distinguishing.\medskip

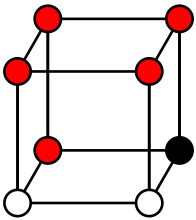
\begin{figure}[htb]
\begin{center}
 \begin{tikzpicture}[scale=1.75]
	 \tikzstyle{edge} = [draw,line width=1pt,-] 
	 \draw[edge] (0,0) -- (0,1) -- (1,1) -- (1,0) -- (0,0);
	 \draw[edge] (0,0) -- (0.23, 0.4) -- (0.23,1.4) -- (0,1) -- (0,0);
	 \draw[edge] (1,0) -- (1.23,0.4) -- (1.23,1.4) -- (1,1) -- (1,0);
	 \draw[edge] (0.23, 0.4) -- (1.23,0.4) -- (1.23,1.4) -- (0.23,1.4) -- (0.23, 0.4);
	  \draw[fill=white!100,line width=1] (0,0) circle (.1); 
	  \draw[fill=red!100,line width=1] (0,1) circle (.1); 
	  \draw[fill=white!100,line width=1] (1,0) circle (.1); 
	  \draw[fill=red!100,line width=1] (1,1) circle (.1); 
	  \draw[fill=red!100,line width=1] (0.23, 0.4) circle (.1); 
	  \draw[fill=red!100,line width=1] (0.23,1.4) circle (.1); 
	  \draw[fill=black!100,line width=1] (1.23,0.4) circle (.1); 
	  \draw[fill=red!100,line width=1] (1.23,1.4) circle (.1);
 \end{tikzpicture}

 \end{center}
 \caption{A 3-distinguishing coloring of $Q_3$}\label{fig:Q3a}
 
 \end{figure}

\begin{example} \label{ex:Q3b} We know that $\det(Q_3)=3$.  Two non-isomorphic minimum determining sets for $Q_3$ are $S_1=\{000,101,110\}$ and $S_2=\{000,010,110\}$.  One way to see that they are non-isomorphic is to look at the the size of their set stabilizers.  We can show that $\setstab(S_1) = \symm(S_1)$ has order 6, while $\setstab(S_2) = \symm(\{000,110\})$ has order 2.  If we use $S_3$ to achieve a distinguishing coloring for $Q_3$, we  need 3 colors to distinguish $S_1$ and therefore 4 colors to distinguish $Q_3$.  However, we can distinguish $S_2$ with 2 colors, allowing us to 3-distinguish $Q_3$.  That is the best we can do, so $\fdist(Q_3)=3$. See Figure~\ref{fig:Q3a} for a 3-distinguishing coloring of $Q_3$ achieved by 2-distinguishing $S_2$.\end{example}

\begin{example} \label{ex:Q8a} We know by  \cite{BoCo2004} that $\dist(Q_8)=2$, by \cite{Bo2013b} that $\rho^2(Q_8) = \urho(Q_8)=5$, by  \cite{Bo2013a}  that $\det(Q_8) = 4$, by \cite{Bo2013b} that all minimum determining sets for $Q_8$ are isomorphic, and each has set stabilizer isomorphic to $\symm(\{1,2,3,4\})$. This implies that every minimum determining set for $Q_8$ requires $4$ colors to distinguish.  Thus $\fdist(Q_8)=5$.\end{example}

Though in the previous two examples the difference between $\rho^u(G)$ and $\rho^\ell(G)$ were small or nonexistent, this difference can be as large as we want.  The next example shows that $\urho(G)$ and $\lrho(G)$, as well as $\fdist(G)$ and $\dist(G)$, can be arbitrarily large multiples of each other.

\begin{example} Let $H$ be an asymmetric graph on $m$ vertices. By \emph{asymmetric} we mean that $\aut(H)$ is trivial.  In \cite{BoIm2017}, to prove that the cost of 2-distinguishing can be an arbitrarily large multiple of the determining number, Boutin and Imrich show that \begin{itemize} \item $\dist(K_{2^m}\Box H)=2,$  \item  $\rho(K_{2^m}\Box H)=m2^{m{-}1}= \frac{|V(K_{2^m}\Box H)|}{2}, \text{ and }$ \item  $\det(K_{2^m}\Box H) = 2^m{-}1.$\end{itemize}  In particular, $\rho(K_{2^m}\Box H)> \frac{m}{2} \cdot\det(K_{2^m}\Box H)$.  In our new notation this says that $\urho(K_{2^m}\Box H)> \frac{m}{2} \cdot\lrho(K_{2^m}\Box H)$.\medskip

In the following theorem, the overall result from \cite{BoIm2017}  is restated using the parameters from this paper.

\begin{thm} \rm   \cite{BoIm2017} To any positive integer $n$ there are infinitely many graphs $G$ for which $\urho(G) >n\cdot \lrho(G)$.\end{thm} 

Next we look at the number of colors necessary to achieve this cost reduction.  But first, let us look at some details of the Cartesian product, $K_{2^m}\Box H$, and its automorphisms. Since $K_{2^m}$ and $H$ are relatively prime, and $H$ is asymmetric, each $\alpha \in \aut(K_{2^m}\Box H)$ has  the form $\alpha(z, h) = (\alpha(z), h)$ for some $\alpha \in \aut(K_{2^m}).$   For each $z\in K_{2^m}$, the set $H^z = \{(z, h) \ | \ h \in H\}$ is called an \emph{$H$-fiber}.  There are $|V(K_{2^m})|=2^m$ $H$-fibers in $K_{2^m}\Box H$.  Further $\alpha \in \aut(K_{2^m} \Box H)$ permutes the $H$-fibers, taking $H^z = \{ (z,h) \ | \ h \in H\}$ to $H^{\alpha(z)} = \{(\alpha(z), h) \ | \ h\in H\}$.  Similarly, the sets  $K_{2^m}^{h} = \{ (z,h) \ | \ z \in K_{2^m}\}$ are called \emph{$K_{2^m}$-fibers}, there are $|V(H)|=m$ of these, and each is preserved setwise by automorphisms of $K_{2^m}\Box H$.

\begin{thm}\rm Let $H$ be an asymmetric graph on $m$ vertices.  Then $$\fdist(K_{2^m}\Box H)=\left\lceil \frac{2^m{-}1}{m}\right\rceil+1.$$\end{thm}

\begin{proof} To find the frugal distinguishing number of $K_{2^m}\Box H$, we start with a minimum determining set and distinguish it. A minimum determining set for $K_{2^m}\Box H$ contains precisely one vertex from all but one of the $2^m$ $H${-}fibers.  This requires $2^m{-}1$ vertices.  However, it doesn't matter in which $K_{2^m}${-}fibers they reside, as long as they reside in distinct $H${-}fibers. Let $S$ be a minimum determining set for $K_{2^m}\Box H$.  For each $h\in V(H)$, let $S^h=S\cap K_{2^m}^h$.  Then $\{S^h\}_{h\in H}$ is a partition of $S$. For each $h$, $\setstab(S^h)=\symm(S^h)$ and is contained in $\setstab(S)$.   Further, the action of $\symm(S^h)$ on $S^h$ requires $|S^h|$ colors to distinguish.   Since we wish to minimize $\max\{S^h\}_{h\in H}$. the number colors necessary to distinguish $S$, we wish to minimize $\max\{|S^h|\}_{h\in H}$.   We achieve this by distributing vertices of $S$ as evenly as possible across all $K_{2^m}$-fibers.  This yields $\max\{|S^h|\}_{h\in H} = \left\lceil \frac{2^m{-}1}{m}\right\rceil$. Since the $K_{2^m}$-fibers are preserved setwise by every $\alpha\in \aut(K_{2^m}\Box H)$, each $S^h\subset K_{2^m}^h$ is preserved setwise by each $\alpha\in \setstab(S)$.  Thus, we can use this same set of $ \left\lceil \frac{2^m{-}1}{m}\right\rceil$ colors on each $S^h$ and we get $\fdist(K_{2^m}\Box H)=\left\lceil \frac{2^m{-}1}{m}\right\rceil+1$.\end{proof}
 
\begin{corollary} \rm  To any positive integer $n$ there are infinitely many graphs $G$ for which $\fdist(G)>n\cdot \dist(G)$.\end{corollary}

\begin{proof}  Since $\dist(K_{2^m}\Box H)=2$ and $\fdist(K_{2^m}\Box H)=\left\lceil \frac{2^m{-}1}{m}\right\rceil+1$, the result is immediate.\end{proof}

\end{example}

\section{Book Graphs}\label{sec:books}

A book graph $B_{m,n}$ consists of $n$ copies of the cycle $C_m$ identified along a single edge.  The individual $m${-}cycles are the {\it pages} of the graph, and the common edge is the {\it spine}. Identify the vertices belonging to the common edge by $v_0,v_{m{-}1}$, and the remaining vertices (in order) on the $i\th$ copy of $C_m$ by $v_{1,i},\ldots, v_{m{-}2,i}$.  Ignoring the vertices on the spine because they are common to all pages, we consider each page as a path on $m{-}2$ vertices. Thus in the following, $v_0$ and $v_{m-1}$ are \emph{spinal-vertices}, the graph induced by the non-spinal-vertices on a page is a \emph{path}, and vertices that are not $v_0$ or $v_{m-1}$ are \emph{path-vertices}.\medskip 

Note that $B_{m,1}=C_m$, which is already well{-}understood, and so $n=1$ is not addressed in this paper.  For $n\geq 2$, the symmetries of $B_{m,n}$ are generated by the $n!$ permutation of the paths, and the nontrivial reflection through $v_0v_{m{-}1}$. In the case $m=3$, the paths are trivial paths in that each consists of a single vertex. The symmetry parameter results for $B_{3,n}$ are different, and arguably less interesting, than those for $B_{m,n}$ with $m\geq 4$.  In the following, we first cover the case $m\geq 4$, and then for completeness sake, the case $m=3$.

\subsection{\boldmath Parameters for $B_{m,n}$ with $m\geq 4$} \label{sec:geq4} 

\begin{thm} \label{thm:det} \rm  If $n\geq 2$ and $m\geq 4$, then $\det(B_{m,n}) = n{-}1$. \end{thm}

\begin{proof}A determining set $S$ for $B_{m,n}$ must contain one vertex from all but one of its paths; otherwise there are at least two paths without vertices in $S$ which can be permuted nontrivially while still fixing $S$ pointwise. Thus $|S|\geq n{-}1$. If $m$ is odd, a determining set must also contain at least one vertex that is not $v_{\frac{m{-}1}{2}, i}$ for any $i$; otherwise the nontrivial reflection through $v_0v_{m{-}1}$ fixes $S$ pointwise. Since $m\geq 4$, we have $1\ne \frac{m{-}1}{2}$, and so $S=\{v_{1,1}, v_{1,2}, \ldots, v_{1,n{-}1}\}$ is a determining set of size $n{-}1$.  Thus $\det(B_{m,n})=n{-}1$.\end{proof} 

In \cite{LaBh2009}, Lal and Bhattacharjya prove the following.

\begin{thm} \label{thm:LB} \rm \cite{LaBh2009} If $n\geq 2$, $m\geq 4$, and $(k{-}1)^{m{-}2}+1\leq n \leq k^{m{-}2}$, then $\dist(B_{m,n})= k$.\end{thm}

Their proof is useful in understanding how to distinguish book graphs.  Here is the gist of it. To $d$-distinguish $B_{m,n}$ it is necessary to assign a distinct $d$-coloring to each of its paths. We call the colorings assigned to these paths \emph{path-colorings}.  The distinctness of the path-colorings breaks the automorphisms that permute paths.  Further, if there is some path-coloring $c_1c_2\ldots c_{m{-}3}c_{m{-}2}$ whose \emph{reversal}  $c_{m{-}2}c_{m{-}3}\ldots c_{2}c_{1}$ is not a path-coloring, these colorings break the reflection as well.  In this case, we can color $v_0$ and $v_{m{-}1}$ the same color if we wish.  On the other hand, if for every path-coloring the reverse path-coloring exists on some path, then we must color $v_0$ and $v_{m{-}1}$ with distinct colors to break the reflection.   We now have the background to study the paint cost of $B_{m,n}$.  Our first theorem applies to $B_{m,n}$ when $n$ is \lq large enough.'

\begin{thm}\label{thm:biggern}\rm If $n\geq 2$, $m\geq 4$, and $d\geq \dist(B_{m,n})$,  then 
$$\rho^d(B_{m,n}) = \left\{ \begin{array}{lll} (m{-}2)(n-d^{m-3}) & \text{ if } &  d^{m{-}2} {-} (d{-}1)^{m{-}2}< n < d^{m-2} \\

(m{-}2)(n-d^{m-3}){+}1 & \text{ if } & n \in \{d^{m{-}2} {-} (d{-}1)^{m{-}2}, d^{m-2}\} \end{array}\right. .$$\end{thm}

\begin{proof} To find $\rho^d(B_{m,n})$ we find the size of the largest possible color class in any $d${-}distinguishing coloring of $B_{m,n}$, and subtract that size from $|V(B_{m,n})|$. In the following we will use \emph{red} as the color whose use we maximize.  Note that the work here is valid for any $d\geq \dist(B_{m,n})$, not just for $\dist(B_{m,n})$ itself.\medskip

Note that there are precisely $d^{m{-}2}{-}(d{-}1)^{m{-}2}$ distinct $d$-path-colorings that use red (all $d$-path-colorings minus the $d$-path-colorings that do not use red).  Thus if $n\geq d^{m{-}2}{-}(d{-}1)^{m{-}2}$, then we can create $n$ distinct path-colorings using the maximum number of red path-vertices. Note that there are ${m{-}2\choose i} {(d{-}1)}^{i}$ path-colorings using precisely $i$ non-red vertices (and therefore precisely $m{-}2{-}i$ red vertices). Thus the total number of red vertices in all path-colorings that use red is $\sum_{i=0}^{m{-}2} (m{-}2{-}i) {m{-}2\choose i} (d{-}1)^i$.  A little binomial algebra shows this total is $(m{-}2)d^{m{-}3}$.\medskip

When $d^{m{-}2} {-} (d{-}1)^{m{-}2}< n < d^{m-2}$, we may use as one of our path-colorings, one that does not use red and whose reversal is not also used as a path-coloring.  In this case, the path-colorings break the reflection, so we can color both spinal-vertices red as well.  This uses the maximum of red vertices, precisely $2{+}(m{-}2)d^{m{-}3}$, in a $d$-distinguishing coloring of $B_{m,n}$.  However, if $n=d^{m{-}2} {-} (d{-}1)^{m{-}2}$, then to maximize red vertices we must use all the $d$-path-colorings that use red, and this set of colorings is preserved under the reflection.  We must  break the reflection by  coloring exactly one spinal-vertex red.   Similarly, if $n=d^{m{-}2}$, we must use all possible $d$-path-colorings, and this set is also preserved under the reflection.  Thus for these two cases, the $d$-distinguishing coloring uses the maximum number of red vertices, precisely $1 {+} (m{-}2)d^{m{-}3}$ .  We subtract this maximum number of red vertices from the total number of vertices, $|V(B_{m,n})|=2+n(m-2)$, and the result follows.\end{proof}

When $n<d^{m{-}2}{-}(d{-}1)^{m{-}2}$, we cannot use all possible $d$-path-colorings that use red, so our summation does not simplify as nicely.  To make it more feasible to state results in this situation, it is useful to have notation for smaller sums.  Let $n^d_j = \sum_{i=0}^j {m{-}2 \choose i} (d{-}1)^i$, which is the number of $d$-path-colorings with at most $j$ non-red vertices per path.  Further, let $N^d_j = \sum_{i=0}^j (m{-}2{-}i){m{-}2 \choose i} (d{-}1)^i$, the number of  red vertices used in $d$-path-coloring $n_j$ paths using at most $j$ non-red vertices per path.\medskip

\begin{thm}\rm \label{thm:smallern} If $d\geq \dist(B_{m,n})$, and $j$ is the largest integer so that $n\leq n^d_j$,  then $|V(B_{m,n})| - N^d_j -1 \leq \rho^d(B_{m,n}) < |V(B_{m,n})| - N^d_{j-1} - 1 $.\end{thm}

\begin{proof}  To be sure that red is used as often as possible, we color all the vertices red on the first path, then all but one vertex red on each of next ${m{-}2\choose 1}(d{-}1)^1$ paths, then all but 2 vertices red on the next ${m{-}2\choose 2}(d{-}1)^2$ paths, and continue in this fashion until we color all $n$ paths.  Since $n_{j-1}^d < n \leq n_j^d$, by the definitions of $N^d_j$, we can bound the maximum number of red path-vertices between $N^d_{j-1}$ and $N^d_j$.  However, we must also adjust for red spinal-vertices.  As argued in the proof of Theorem~\ref{thm:biggern}, if $n=n_j^d$ we use every path-coloring containing at most $j$ non-red vertices, and this set of path-colorings is preserved under the reflection automorphism. In this situation, we color only one of the spinal-vertices red, which breaks the reflection.   Thus $R^d$, the maximum number of red vertices used in a $d$-distinguishing coloring of $B_{m,n}$, fits the inequality $1{+}N^d_{j-1} < R^d \leq 1{+}N^d_j$.  The result follows. \end{proof}

The following two corollaries reframe the above results in terms of $\dist(B_{m,n})$, which allows us to claim bounds on $\urho(B_{m,n})$.  The first corollary gives broad bounds, applicable to all $n$, and is not too hard to state.  The second gives more refined bounds, but requires the notation of $n_j^d$ and $N_j^d$. 

\begin{corollary}  \rm If $n\geq 2$, $m\geq 4$, and  $k=\dist(B_{m,n})$,  then $$(m{-}2)(n- k^{m{-}3})+1 \leq \rho^u(B_{m,n}) <(m{-}2)(n- (k-1)^{m{-}3})+1.$$\end{corollary}\medskip

\begin{corollary}\rm If $n\geq 2$, $m\geq 4$, $k=\dist(B_{m,n})$, and $j$ is the largest integer so that $n\leq n^d_j$,  then $$|V(B_{m,n})| {-} N^k_j {-}1 \leq \rho^u(B_{m,n}) < |V(B_{m,n})| {-} N^k_{j-1} {-} 1.$$\end{corollary}

Now let us find $\fdist(B_{m,n})$. 

\begin{thm}\rm If $n\geq 2$ and $m\geq 4$, then $\fdist(B_{m,n}) = 2+\left\lfloor \frac{n{-}1}{m{-}2} \right\rfloor$. \end{thm}

\begin{proof} We argued in the proof of Theorem~\ref{thm:det}  that for $m\geq 4$, $S$ is a minimum determining set for $B_{m,n}$ if and only if it contains exactly one path-vertex from all but one of the paths, and at least one of the path-vertices of $S$ is not the center vertex of its associated path. Choose such a determining set $S$.  Partition $S$ into subsets containing vertices in the same position on their respective paths. That is, for each $j\in \{1,\ldots, m{-}2\}$, let $S_j=\{ v_{j,i} \ | \  v_{j,i} \in S\}$.  By the definition of $S_j$, all possible permutations of paths associated with vertices of $S_j$ are contained in $\setstab(S)$.  Thus the set $S_j$ requires $|S_j|$ colors to be distinguished, and $S$ requires at least $\max\{|S_j|\}_{j=1}^{m{-}2}$ colors.  Since  $\aut(B_{m,n})$, and therefore $\setstab(S)$, is composed of path permutations and (possibly) a  reflection, $S$ requires precisely $\max\{|S_j|\}_{j=1}^{m{-}2}$ colors to distinguish as long as the reflection does not preserve $S$ setwise.  In the following, we choose a minimum determining set $S$ that is not preserved by the reflection.\medskip

Let $q=\lfloor \frac{n{-}1}{m{-}2} \rfloor$.  Then $n{-}1=q(m{-}2){+}r$ with $0\leq r< m{-}2$.  Distribute the $m{-}2$ path positions over $n{-}1$ paths as evenly as possible, achieving minimum determining set $S$.  Partition $S$ into $\{S_j\}_{j=1}^{m{-}2}$ as described above.  If $r=0$,  then for each $j$, $|S_j|=q$.  If $r>0$, then  without loss of generality, or by reordering the paths, for $j\leq r$, $|S_j|=q{+}1$ while for $j>r$, $|S_j|=q$.\medskip

We now show that $S$ can be $(q{+}1)$-distinguished.  Suppose that $r>0$.  Then $|S_1|=q{+}1$, while $|S_{m{-}2}|=q$, and so the reflection automorphism does not preserve $S$.  Thus, automorphisms in $\setstab(S)$ are composed strictly of permutations of the paths associated with the same $S_j$.  These permutations can be distinguished using $\max\{|S_j|\}_{j=1}^{m{-}2} = q{+}1$ colors. Thus, we have a $(q{+}1)$-distinguishing coloring of $S$.\medskip

Suppose that $r=0$.  In this case we can break all page permutations by using $q$ colors to distinguish each of the $S_j$, but then the reflection automorphism preserves the color classes.  Thus we either need to use $q{+}1$ colors to distinguish $S$ so that say $S_1$ is distinguished using $\{1,\ldots, q\}$ and $S_{m-2}$ is distinguished using $\{2, \ldots, q+1\}$, or we need to modify our partition so that one part has cardinality $q{+}1$, another $q{-}1$, and the remainder $q$.  Either solution requires $q{+}1$ colors to distinguish $S$, and therefore $q{+}2$ colors to achieve the lower paint cost $\det(B_{m,n})$.\end{proof}

\begin{example}  \rm If $m=8$ and $n=473$, since $472=78\cdot 6 + 4$, we can $80${-}distinguish $B_{8,473}$ at a cost of $\rho^\ell(B_{8,473})=472$.  Alternatively, since $\dist(B_{8,473})=3$, we can 3{-}distinguish $B_{8,473}$ at a cost of $\rho^u(B_{8,473})=1573$. 

If $n=703$, since $700=116\cdot 6 + 7$, we can $118${-}distinguish $B_{8,703}$ at a cost of $\rho^\ell(B_{8,703})=702$.  Alternatively, since $\dist(B_{8,703})=3$ if we 3{-}distinguish $B_{8,703}$, the number of red vertices maxes out at 1458, the number of vertices is 4220, and thus $\rho^u(B_{8,703}) = 2762$.   \end{example}
 
 \subsection{\boldmath Parameters for $B_{3,n}$}
 
 Let $m=3$.  In this case, the \lq\lq paths" defined by deleting vertices $v_0$ and $v_{m-1}$ are isolated vertices.  However, we continue the terminology of paths, while understanding that here they are paths only trivially. 
 
 \begin{thm}\rm Let $n\geq 2$.  Then
 
 \begin{enumerate}
 
 \item[(i)] $\det(B_{3,n}) = n$,
  
 \item[(iii)] $\rho^u(B_{3,n})=n$, and
 
 \item[(iv)] $\fdist(B_{3,n})= n$. 
 
 \end{enumerate} \end{thm} 
 
 \begin{proof} As argued in the proof of Theorem~\ref{thm:det}, a determining set for $B_{3,m}$ must contain one vertex from all but one of the paths.  Since $m=3$, all path-vertices are of the form $v_{i,1}$, and are fixed by the nontrivial reflection.  Thus, it is not sufficient for a minimum determining set $S$ to have only path-vertices; it is necessary to add one of the spinal-vertices to account for the reflection automorphism.  Thus $S=\{v_0,v_{1,1}, v_{2,1}, \ldots, v_{n{-}1, 1}\}$ is a minimum determining set for $B_{3,n}$ of size $n$.\medskip
 
 Note that all minimum determining sets for $B_{3,n}$ are isomorphic. Each must contain $n-1$ vertices of the form $v_{1,i}$ and one of the spinal-vertices. Let $S$ be such a determining set. Since $\setstab(S) = \symm(\{v_{1,i}\}_{i=1}^{n-1}$, $S$ requires $n{-}1$ colors to distinguish its path-vertices, while the unique spinal-vertex in $S$ can be colored with the same color as one of its path-vertices.  Thus every minimum determining set requires $n{-}1$ colors to distinguish and therefore $\fdist(B_{3,n})=n$.  By  Lal and Bhattacharjya~\cite{LaBh2009}, $\dist(B_{3,n})=n$ as well.  Thus $\fdist(B_{3,n}) = \dist(B_{3,n})$ and therefore $\rho^\ell(B_{3,n}) = \rho^u(B_{3,n})=n$.
 \end{proof}
 
 \section{Future Work}\label{sec:future}
 
 Since these parameters are new, the area is wide open for investigation.  However, the following are some specific open problems.
 
 \begin{problem} Revisit families of 2-distinguishable graphs (some example given in the introduction) to investigate  their upper paint costs and their frugal distinguishing numbers.  How does these differ within or between graph families? \end{problem}
 
 \begin{problem}  Classify graphs for which $\dist(G)=\fdist(G)$. \end{problem}
 
 \begin{problem}  Classify graphs for which $\dist(G)$ and $\fdist(G)$ differ by at most a fixed constant, or a fixed multiplier. \end{problem}
 
\bibliographystyle{plain}
\bibliography{FrugalDist}

\begin{thebibliography}{10}

\bibitem{Al2005}
Michael~O. Albertson.
\newblock Distinguishing {C}artesian powers of graphs.
\newblock {\em Electron. J. Combin.}, 12:Note 17 (electronic), 2005.

\bibitem{AlBo2007}
Michael~O. Albertson and Debra~L. Boutin.
\newblock Using determining sets to distinguish {K}neser graphs.
\newblock {\em Electron. J. Combin.}, 14(1):Research Paper 20 (electronic),
  2007.

\bibitem{AlCo1996}
Michael~O. Albertson and Karen~L. Collins.
\newblock Symmetry breaking in graphs.
\newblock {\em Electron. J. Combin.}, 3(1):Research Paper 18, approx. 17, 1996.

\bibitem{AlSo2021a}
Saeid Alikhani and Samaneh Soltani.
\newblock The cost number and the determining number of a graph, 2021.
\newblock arXiv:1710.07527.

\bibitem{Ba1977}
L.~Babai.
\newblock Asymmetric trees with two prescribed degrees.
\newblock {\em Acta Math. Acad. Sci. Hungar.}, 29(1-2):193--200, 1977.

\bibitem{BoCo2004}
Bill Bogstad and Lenore~J. Cowen.
\newblock The distinguishing number of the hypercube.
\newblock {\em Discrete Math.}, 283(1-3):29--35, 2004.

\bibitem{Bo2020}
Debra Boutin.
\newblock The cost of 2-distinguishing hypercubes, 2020.
\newblock arXiv:2007.15948.

\bibitem{BoIm2017}
Debra Boutin and Wilfried Imrich.
\newblock The cost of distinguishing graphs.
\newblock In {\em Groups, graphs and random walks}, volume 436 of {\em London
  Math. Soc. Lecture Note Ser.}, pages 104--119. Cambridge Univ. Press,
  Cambridge, 2017.

\bibitem{Bo2008}
Debra~L. Boutin.
\newblock Small label classes in 2-distinguishing labelings.
\newblock {\em Ars Math. Contemp.}, 1(2):154--164, 2008.

\bibitem{Bo2013a}
Debra~L. Boutin.
\newblock The cost of 2-distinguishing {C}artesian powers.
\newblock {\em Electron. J. Combin.}, 20(1):Paper 74, 13, 2013.

\bibitem{Bo2013b}
Debra~L. Boutin.
\newblock The cost of 2-distinguishing selected {K}neser graphs and hypercubes.
\newblock {\em J. Combin. Math. Combin. Comput.}, 85:161--171, 2013.

\bibitem{ImKl2006}
Wilfried Imrich and Sandi Klav{\v{z}}ar.
\newblock Distinguishing {C}artesian powers of graphs.
\newblock {\em J. Graph Theory}, 53(3):250--260, 2006.

\bibitem{ImKlTr2007}
Wilfried Imrich, Sandi Klav{\v{z}}ar, and Vladimir Trofimov.
\newblock Distinguishing infinite graphs.
\newblock {\em Electron. J. Combin.}, 14(1):Research Paper 36, 12 pp.
  (electronic), 2007.

\bibitem{ILTW2020}
Wilfried Imrich, Thomas Lachmann, Thomas~W. Tucker, and Gundelinde~M. Wiegel.
\newblock Asymmetrizing cost and density of vertex-transitive cubic graphs.
\newblock Personal Copy from Author, 2020.

\bibitem{KlZh2007}
Sandi Klav{\v{z}}ar and Xuding Zhu.
\newblock Cartesian powers of graphs can be distinguished by two labels.
\newblock {\em European J. Combin.}, 28(1):303--310, 2007.

\bibitem{LaBh2009}
A.~K. Lal and B.~Bhattacharjya.
\newblock Breaking the symmetries of the book graph and the generalized
  {P}etersen graph.
\newblock {\em SIAM J. Discrete Math.}, 23(3):1200--1216, 2009.

\bibitem{SmTuWa2012}
Simon~M. Smith, Thomas~W. Tucker, and Mark~E. Watkins.
\newblock Distinguishability of infinite groups and graphs.
\newblock {\em Electron. J. Combin.}, 19(2):Paper 27, 10, 2012.

\bibitem{WaZh2007}
Mark~E. Watkins and Xiangqian Zhou.
\newblock Distinguishability of locally finite trees.
\newblock {\em Electron. J. Combin.}, 14(1):Research Paper 29, 10 pp.
  (electronic), 2007.

\end{thebibliography}

\end{document}